\documentclass[11pt]{amsart}
\usepackage{amssymb,amsmath,latexsym,amscd,amsfonts}
\usepackage{graphics}
\usepackage{epsfig}
\usepackage{psfrag}
\def\ignore #1 {}

\newtheorem{lem}{Lemma}
\newtheorem{thm}{Theorem}

\newtheorem{cor}{Corollary}
\def\hpic #1 #2 {\mbox{$\begin{array}[c]{l} \epsfig{file=#1,height=#2} \end{arr\
ay}$}}
\def\vpic #1 #2 {\mbox{$\begin{array}[c]{l} \epsfig{file=#1,width=#2} \end{arra\
y}$}}

\def\ga{\gamma}
\def\al{\alpha}
\def\pr{\mbox{Prob}}

\def\F{\mbox{{\bf F}}}
\def\N{\mbox{{\bf N}}}

\def\R{\mbox{{\bf R}}}

\def\Z{\mbox{{$\mathbb{Z}$}}}

\def\a{\mbox{${\alpha}$}}

\renewcommand{\a}{\alpha}
\newcommand{\g}{\gamma}
\newcommand{\e}{\epsilon}
\renewcommand{\b}{\beta}
\renewcommand{\d}{\delta}

\begin{document}
\title{Random Multiplication Approaches Uniform Measure in Finite Groups}
\author{Aaron Abrams, Henry Landau, Zeph Landau,\\
James Pommersheim, Eric Zaslow}

\begin{abstract}
In order to study how well a finite group might be generated by repeated random 
multiplications, P.~Diaconis suggested the following urn model.  An urn contains 
some balls labeled by elements which generate a group $G$.  Two are drawn at 
random with replacement and a ball labeled with the group product (in the order 
they were picked) is added to the urn.  We give a proof of his conjecture that the 
limiting fraction of balls labeled by each group element almost surely approaches 
$\frac{1}{|G|}$.
\end{abstract}


\maketitle
\section{Introduction}
In\footnotetext{AA supported in
part by NSF grant DMS-0089927.
The work of EZ has been supported in part by an
Alfred P. Sloan Foundation fellowship.  EZ also thanks
the School of Mathematics, Statistics, and Computer
Science at Victoria University, Wellington, where some
of this work took place.} order to study how well a finite group might be generated by repeated random multiplications, 
P.~Diaconis suggested the following urn model.  An urn contains some balls labeled by elements which
generate a group $G$.  Two are drawn at random with replacement and a ball labeled with the group product
(in the order they were picked) is added to the urn. He conjectured that the limiting fraction of balls labeled by 
each group element approaches $\frac{1}{|G|}$ with probability 1.

This problem arose from work of Diaconis and S.~Rees who were studying a group theoretic
algorithm called MeatAxe.  (For further reading see~\cite{holt-rees}.)  This is a widely used 
tool for decomposing representations of a 
finite group $G$ over a finite field $\F$.  To begin the MeatAxe, a random element of the group
algebra $\F G$ must be chosen.  In practice, this is done by taking a sum of a few products of 
generators such as $B+AB+BA+BAB^2+A+A$ and ``hoping for the best.''

Diaconis and Rees began the study of more careful algorithms which would provably converge
to a random element of $\F G$.  One proposal was this:  from $x\in \F G$, go to $sx$, $s^{-1}x$,
or $x+as$, where $s$ is uniformly chosen from a generating set of $G$
and $a$ is uniformly chosen from $\F$.  The problem studied here arose as a sub-problem in
analyzing this algorithm; it turns out to be challenging even for a group with two elements.

Note that, if $G= \Z /2\Z =\{0,1\}$ and the urn contains a large fraction of $0$'s, then the probability of adding 
another $0$ is also large, and it may seem that the preponderance of $0$'s continues.  On the other hand, 
a simple computation shows that the expectation is that the fraction of $1$'s moves toward 1/2.

In this paper we use an elementary method to prove the conjecture of Diaconis.
D.~Siegmund and B.~Yakir~\cite{siegmund} have obtained these results independently
as an application of an almost supermartingale convergence theorem, and another proof based 
on large deviations has been given by A.~Shwartz and A.~Weiss~\cite{shwartz-weiss}.
We remark that the rate of convergence of this procedure
is open at this writing.

The work of EZ has been supported in part by an
Alfred P. Sloan Foundation fellowship.  EZ also thanks
the School of Mathematics, Statistics, and Computer
Science at Victoria University, Wellington, where some
of this work took place. 
The authors thank P.~Diaconis
for introducing us to this problem and for helpful conversations along the way.  We
also thank Julie Landau whose support was essential to this work.


\section{Notation and Outline}

Let $G=\{g_1,\ldots,g_d\}$ be a finite group.  The state $s$ of an
urn is described by the number of balls with each
of the $d$ possible labels; thus we write 
\begin{equation*}
s=(n_1,\ldots,n_d)
\end{equation*}
where $n_i=n_{g_i}$ is the number of balls labeled $g_i$.  
We measure time by the total number of balls:
\begin{equation*}
t=\sum_{i=1}^d n_i(t);
\end{equation*}
in particular the starting time of the process is a positive integer
determined by the initial configuration.  To emphasize that the state
evolves we will write $s=s(t)$ and $n_i=n_i(t)$.  Let 
\begin{equation*}
p_i(s(t))=p_{g_i}(s(t)):=\frac{n_i(t)} t 
\end{equation*}
be the fraction of balls 
labeled $g_i$ in state $s(t)$; we will also refer to $p_i(s)$ as the 
\emph{density} of such balls.  Clearly, $\sum_i p_i(s)=1$.  We often 
omit the explicit dependence on $s$ and $t$ and write $p_i$ for 
$p_i(s(t))$ and $n_i$ for $n_i(t)$.

Our main result is the following.

{\bf Corollary \ref{main}.}  
\emph{Assume the urn initially contains a set of balls whose labels generate
the finite group $G$.  Then with probability one, for each $i$,
$$ p_i \longrightarrow 1/d$$
as $t \to \infty$.
}

In other words, in the limit, the urn contains an equal fraction of all
group elements.

We now describe some of the notation and strategy used in the proof.
Our process can be described 
as a random walk on the non-negative integer lattice in $\R^d$, 
with the (state-dependent) probability of adding one to the $j$th 
coordinate given by
\begin{equation} \label{e:prob}
\pi_j(s) :=\sum _{g\in G} p_g(s) p_{g^{-1}g_j}(s).
\end{equation}
That is, at state $s=(n_1,\ldots,n_d)$ we have
\begin{equation}\label{e:probrule}
n_j(t+1)=
\begin{cases}
n_j(t) +1 & \text{w.p. } \pi_j(s(t)) \\
n_j(t) & \text{otherwise}.
\end{cases}
\end{equation}
Observe, of course, that the functions $n_j$ do not evolve 
independently, as they are constrained by $\sum n_i(t)=t$.

Now, let $0 \leq \a \leq 1/d$, and define 
\begin{equation*} \label{e:sig}
\Sigma_{\al}=\left\{ (x_1, \dots x_d)\in \R^d \mid
x_j \geq \al (x_1+\cdots+x_d) \text{ for each } 1\leq j \leq d\right\};
\end{equation*}
thus $s(t)\in \Sigma_\a$ if the balls of each label have density 
at least $\a$.

We will deduce our main result by showing that as $t$ grows, 
$\pr\{s(t)\in \Sigma_\b\}$ 
almost surely approaches 1 for each $\b<1/d$.
The argument has two main ingredients which are contained in Lemma \ref{l:2} and 
Theorem \ref{l:3}.  Lemma \ref{l:2} shows that with probability exponentially close to 1, 
if the distribution of labels is not uniform then the lowest 
density in the urn increases by a constant factor after evolving 
for a fixed fraction of the elapsed time.  By iterating, we can bring the lowest density
arbitrarily close to $1/d$.  To make this result effective, however, this lowest density 
must be nonzero, i.e., every element of $G$ must have a representative in the urn.  This is clearly true for 
$G= \Z/2\Z$, so Lemma \ref{crux} is sufficient for this case.  However we require the
additional arguments of Section \ref{general}, which use the result in the $\Z/2\Z$ case, to 
show that an urn containing only a generating set will eventually (with probability 1) 
contain each element of $G$.

%

\section{Special case:  the urn contains every element}

We begin with a 
bound on the transition probabilities $\pi_j(s)$
for states $s\in \Sigma_\al$.

\begin{lem} \label{l:1}
If $s\in\Sigma_\al$, then for each $j$, the transition
probability satisfies
\begin{equation*}
\pi_j(s) \geq 2 \al - d \al^2.
\end{equation*}
\end{lem}

\begin{proof}
We prove more generally that if $x_i$, $y_i \geq \al$ for $i=1,\ldots,d$
and $\sum x_i = \sum y_i =1$ then 
$\sum x_iy_i \geq 2\al - d\al^2.$  The lemma follows by setting 
$x_i=p_{g_i}$ and $y_i=p_{g_i^{-1}g_j}$ in (\ref{e:prob}).

Writing $x_i=\a+x_i'$ and $y_i=\a+y_i'$, with $x_i'$, $y_i' \geq 0$, $\sum_{i=1}^d x_i'=\sum_{i=1}^d y_i'=1-d\al$, we
have
\begin{equation*}
\sum x_iy_i=\sum [\a^2 + \a(x_i'+y_i') + x_i'y_i']
= d\a^2+\a(2-2d\a)+\sum x_i'y_i'
\geq 2\a -d\a^2.
\end{equation*}
\end{proof}

The main lemma, which follows, says that (for appropriate
choices of $r, \g >1$) if the state $s$ of the urn
is in $\Sigma_\a$ at time $T$, then at
time $rT$ the chance that the state is in 
$\Sigma_{\g\a}$ is exponentially close to 1.  We want to 
iterate this argument until we can conclude with high probability
that the state is in $\Sigma_{\b}$;
for the iteration it is important that $\g$ and the coefficient
of the exponential term do not depend on $\a$.

\begin{lem} \label{l:2}\label{crux}
Let $0< \al < 1/d$, and fix $r$ in the range $1<r<2-\a d$. 
There exist constants $\g > 1$ and $C>1$, depending only on $r$, 
such that for all $T$:
\begin{equation}\label{e:lemma2}
\mbox{if }\ s(T)\in\Sigma_\a\ \mbox{ then }\ \pr\{s(rT)\in\Sigma_{\g\a}\} \geq 
1-C^{-\a T}.
\end{equation}
\end{lem}

\begin{proof}
We examine the walk independently in each coordinate direction.  
The main idea is that in a given coordinate, we can compare the 
process (\ref{e:probrule}) over a certain period of time to a random walk with 
the constant transition probability given by the bound in Lemma \ref{l:1}. 
As this is a simple random walk, it is easy to estimate its behavior. 
The choice of time period is somewhat delicate,
however:  if it's too short, there isn't enough time to progress
toward the mean (with high probability); on the other hand the 
constant lower bound on the transition probability is only valid for a
limited time, since any $p_i$ could eventually become
arbitrarily close to zero.

Suppose the urn is in state 
$s(T)=(n_1(T),\ldots,n_d(T))\in\Sigma_\a$ at some time $T$.  Note first that, for  $T\leq t \leq rT$, the value of 
$p_i(t)$ can never fall below
\begin{equation} \label{e:overr}
n_i(T)/rT \geq \al/r.
\end{equation}  
As this is true for every $i$, it follows from Lemma \ref{l:1} that 
\begin{equation} \label{e:9}
\pi _j (s(t))\geq  2 \left(\frac{\al}{r}\right) - d \left(\frac{\al}{r}\right) ^2
\end{equation}
for $T\leq t\leq rT$.

Now consider a random walk $W_i$ on the integers which begins at time
$T$ at the value $n_i(T)$, and evolves according to the probabilities
\begin{equation*}
W_i(t+1)=
\begin{cases}
W_i(t)+1 & \text{with probability } 2\a /r - d \a^2/r^2 \\
W_i(t) & \text{otherwise.}
\end{cases}
\end{equation*}

By comparing the evolution (\ref{e:probrule}) with $W_i$ we see from (\ref{e:9}) that for $T\leq t \leq rT$ and any $x$,
\begin{equation*}
\pr \{n_i(t) \geq x \} \geq \pr \{ W_i(t) \geq x \}.
\end{equation*}
It follows that for any $\g$,
\begin{eqnarray} \label{e:6}
\pr \{ p_i(rT) \geq \g \a \} & = & \pr \{\frac{n_i(rT)}{rT} \geq  \g \a\} \\ 
& \geq & \pr \{W_i(rT) \geq  \g\a rT \} \nonumber \\
& =       & \pr \{X_i((r-1) T) \geq \g\a rT -n_i(T) \} \nonumber \\
& \geq & \pr \{X_i((r-1) T) \geq (\g r-1)\a T \}, \nonumber
\end{eqnarray}
where $X_i(t):= W_i(T+t)-n_i(T)$ is the (space and time) translate of $W_i$ 
which starts at 0 at time 0.

This last probability is easy to estimate, as $X_i$ is a sum of 
independent identically distributed Bernoulli trials.  Specifically, the random 
variable $X_i(t+1)-X_i(t)$ takes values 0 and 1 and has mean
\begin{equation}
\mu :=  2 \frac{\al}{r} - d \left(\frac{\al}{r}\right) ^2.
\end{equation}
Then since $d\a<2-r$ and $1<r<2$ we have
\begin{equation}\label{e:muoveral}
\frac{\mu}{\al} =\frac{2}{r} - d \frac{\al}{r^2} > \frac{3r-2}{r^2}>1.
\end{equation}
It follows that $\lim_{\ga \rightarrow 1^+}\frac{(\g r-1)\al}{(r-1)\mu}=\frac{\al}{\mu} <1$, hence 
we can pick $\g>1$ depending on $r$ (but not $\al$) so that 
\[ \d:=\frac{(\g r-1)\al}{(r-1)\mu}<1. \]
The Chernoff bound \cite{chernoff} estimates the probability after time $t$
that $X_i$ is above a fraction $\d$ of expected value:
\begin{equation}\label{chernoff}
\pr\{X_i(t)\geq \d\mu t\} \geq 1-\exp(-(1-\d)^2\mu t/2).
\end{equation}
Applying (\ref{chernoff}) at time $t=(r-1)T$ and with the above $\d$, we obtain
\begin{eqnarray*}
\pr \{X_i((r-1)T)\geq (\g r -1)\a T\} & \geq & 1-\exp\left\{-(1-\d)^2\mu (r-1)T/2\right\} \\
&\geq & 1-\exp\left\{-(1-\d)^2 (r-1) \a T/2\right\}\\
&=&1-A^{-\a T}
\end{eqnarray*}
where the second inequality follows from (\ref{e:muoveral}) and where
\begin{equation*}
A=\exp \{\frac 1 2  (1-\d)^2 (r-1)\}>1.
\end{equation*}

In view of (\ref{e:6}), this gives the desired bound for each $X_i$.
Since the identical bound holds for each $i$, we conclude that the chance 
that the state $s(rT)$ is inside $\Sigma_{\g \a}$ is at least
$1- dA^{-\al T}=1-C^{-\a T}$.
\end{proof}

Note once again that Lemma \ref{crux} provides no information if $\a=0$.  
Theorem \ref{t:1} below assumes $\alpha$ is positive, and Theorem \ref{l:3} 
in the next section shows that this assumption is eventually valid when the 
urn begins with a set of generators for $G$.
We need one more preparatory lemma before proving Theorem \ref{t:1}.
  
\begin{lem}\label{arbitrarilymany}
Suppose there is a time at which $n_i\geq 1$ for each $i$.
Then for each $N\in\N$, there is almost surely a time at which
$n_i\geq N$ for each $i$.
\end{lem}

\begin{proof}  We proceed by induction on $N$.
Assume each $n_i(T)\geq N \geq 1$ at some time
$T$.  Thus for every $t>T$ we have $p_i\geq N/t$, so by Lemma \ref{l:1}, 
\begin{equation*}
\pi_j(s)\geq \frac {2N}{t} - \frac {dN^2}{t^2} > N/t
\end{equation*}
for each $j$, where the last inequality follows since necessarily 
$T\geq dN$.  That is, the chance of increasing $n_j$ at time $t$ 
is at least $N/t$ for each $t>T$.  As $\sum (N/t)$ diverges, with
probability one $n_j$ will eventually increase.  This holds for all $j$
so eventually each $n_j$ will be at least $N+1$.
\end{proof}

\begin{thm} \label{t:1}  Suppose the urn is in a state with each $n_i \geq 1$.  Then 
with probability one, for each $i$,
$$p_i \longrightarrow 1/d$$
as $t\to\infty$.
\end{thm}

\begin{cor}\label{c:1}
Let $G=\Z/2\Z=\{0,1\}$.  If the urn initially contains at least one ball labeled 1,
then the densities of the two elements almost surely approach 1/2.
\end{cor}

\begin{proof}[Proof of Corollary \ref{c:1}.]
We are given that $n_1\geq 1$, and after possibly waiting for one step we will have
$n_0\geq 1$ as well.  Thus Theorem \ref{t:1} applies.
\end{proof}

\begin{proof}[Proof of Theorem \ref{t:1}.]
Fix $\b<1/d$ and let $\e>0$.
It suffices to show that $t'$ can be chosen so that 
$\pr\{s(t)\in\Sigma _\b\}>1-\e$ for all $t>t'$.

Choose $\b'$ between $\b$ and $1/d$.  We will show that with
high probability, the state evolves into $\Sigma_{\b'}$ and then 
stays inside $\Sigma_\b$.  For both steps we will apply 
Lemma \ref{crux} iteratively.

Let $r=2-d\b'$.  Note that for any $\a<\b'$, this choice of $r$ satisfies
the hypothesis of Lemma \ref{crux}; thus we have a $C$
and $\g$ depending on $r$ (but not on $\a$) such that for all
$\a<\b'$ and for all $T$, (\ref{e:lemma2}) holds.

For any $\a<\b'$, then, we may iterate Lemma \ref{crux}
with the same value of $r$.  Thus for $s(T)\in\Sigma_\a$, we
have a bound on the chance that $s(r^jT)$ is \emph{not} in 
$\Sigma_{\g^j \a}$:
\begin{equation}\label{iteratederror}
\pr\{s(r^jT)\not\in\Sigma_{\g^j \a} \}\leq  \sum_{i=0}^{j-1} C^{-\a T(\g r)^i}
\end{equation}
as long as $\g^{j-1}\a<\b'$.
As $j$ tends to infinity, the above sum converges; indeed we
can bound the sum independently of $j$ by a function $f$ of $\a T$
which decreases to 0 as $\a T$ tends to infinity.

Choose $N\in\N$ so large that $f(N)<\frac \e 3$.
By assumption, there is a time at which each $n_i \geq 1$; therefore
by Lemma \ref{arbitrarilymany}
we may choose $T$ such that with probability at least $1-\frac \e 3$, 
each $n_i(T)\geq N$. 
Let $\a=N/T$; then by definition
\begin{equation*}
\pr \{s(T)\in \Sigma_\a (T) \} > 1- \frac \e 3,
\end{equation*}
and so by (\ref{iteratederror}) and our choice of $N$,
\begin{equation}\label{insidebetaprime}
\pr \{s(r^kT)\in\Sigma_{\b'}\} > 1 - \frac {2\e} 3
\end{equation}
where $k$ is the smallest positive integer with $\g^k\a>\b'$.

Let $t'=r^kT$.
It remains to show that the probability of escaping $\Sigma_\b$ at any time 
$t>t'$ can be made arbitrarily small.  For this we apply Lemma 
\ref{crux} with the role of $\a$ played by $\b'$ and the role of $r$
played by any $r'>1$ chosen to be less than $\min\{\b'/\b,2-d\b'\}$.  
The lemma provides a $C'>1$ and a $\g'>1$, and as $\Sigma_{\g' \b'}
\subset \Sigma_{\b'}$, we may ignore $\g'$ and write, for all $t$,
\begin{equation*}
\mbox{If }\ s(t) \in \Sigma_{\b'}\ \mbox{ then }\ \pr\{s(r't)\in\Sigma_{\b'}\} \geq 1-(C')^{-\b't}.
\end{equation*}
Iterating as before, it follows that if $\b't$ is sufficiently large, we can ensure that
\begin{equation}\label{leakagebound}
s(t)\in\Sigma_{\b'}\ \mbox{ implies }
\ \pr\{s((r')^kt)\in\Sigma_{\b'}\ \mbox{ for all }k>0\} > 1-\frac \e 3.
\end{equation}
Note that if $s((r')^kt')\in\Sigma_{\b'}$ for all $k\geq 0$ then $s(t)\in\Sigma_{\b}$
for all $t>t'$, because $\b'/r'>\b$.

So, to complete the proof, we must rechoose our initial value of $T$ if necessary
so that $t'=r^kT$ is large enough for (\ref{leakagebound}) to hold with $t=t'$.
Then by (\ref{insidebetaprime}) and (\ref{leakagebound}) we have
\begin{equation*}
\pr\{s(t)\in\Sigma_{\b}\} > 1- \e
\end{equation*}
for all $t>t'$.
\end{proof}


\section{General case:  the urn contains a generating set}\label{general}

The goal of this section is to prove Theorem \ref{l:3}, from which the main
result (Corollary \ref{main}) follows.

\begin{thm} \label{l:3}
Suppose the urn contains a set of balls whose labels generate the group $G$ 
of order $d$.  Eventually, with probability one, the urn will contain a ball labeled
by each group element.
\end{thm}

\begin{cor} \label{main}
Suppose the urn initially contains a set of balls whose labels generate $G$.  
Then with probability one, for each $i$,
$$p_i \longrightarrow 1/d$$
as $t\to\infty$.
\end{cor}

\begin{proof}[Proof of Corollary \ref{main}.]
This follows immediately from Theorems \ref{t:1} and \ref{l:3}.
\end{proof}

To prove the theorem we will need three lemmas,
which are stated next but proved after the theorem.
In outline, the argument begins at a time $T$ (to be chosen big enough) with a set $S$ of
elements, each with density bigger than a small constant.  Lemma \ref{F1} shows that
with high probability, at a later moment all of the elements of the 
subgroup $H$ generated by $S$ will have density bigger than a smaller constant.  
Lemma \ref{F2} produces a later time when an element
outside $H$ will have a comparable density, and Lemma \ref{F3} shows that at that
time, with high probability the density of the elements of $H$ will not have fallen very 
much.  In this way, a larger set with nonzero density is produced.  By iterating we 
conclude that eventually the whole group will have positive density.

For a subset $S\subseteq G$ let $n_S(t)=\sum_{g\in S} n_g(t)$ and 
$p_S(t)=\sum_{g\in S} p_g(t)=n_S(t)/t$.  Also let $G\setminus S$ denote
the complement of $S$ in $G$, and let $\langle S \rangle$
denote the subgroup of $G$ generated by $S$.

\begin{lem} \label{F1}
Let $0\leq\nu\leq 1$.  There exists a constant $c>1$, depending on $\nu$,
such that for all $T$, if $p_g(T)\geq \nu$ for all $g$ in a subset $S$ of $G$,
then with probability at least $1-c^{-T}$,
$$p_g(T_1) \geq 16\left(\frac{\nu}{16}\right)^{2^d} \quad \mbox{ for all }
g\in\langle S \rangle,$$
where $T_1=2^dT$.
\end{lem}

\begin{lem}\label{F2} 
Let $H\ne G$ be a subgroup of $G$, and let $T$ be any time.  With probability $1$ 
there exists a time $t>T$ for which $p_{G\setminus H}(t)\geq \frac{1}{4}$.
\end{lem}

\begin{lem} \label{F3}
For $\nu$ sufficiently small, there exists a constant $C$ depending on $\nu$
such that for all $T$ and $T'$, if $p_h(T)\geq\nu$ for all $h$ in a subgroup $H$ and
$p_{G\setminus H}(t) < d\nu$ for all $T\leq t < T'$, then with probability at least
$1-C^{-T}$,
$$p_h(T')\geq \frac{2\nu}{3} \quad \mbox{ for all }h\in H.$$
\end{lem}

\begin{proof}[Proof of Theorem \ref{l:3}.]
Let $\nu_0<1/d$ be small enough that Lemma \ref{F3} holds and also that 
$\nu_0'=16 \left(\frac{\nu_0}{16}\right)^{2^d}<\frac{1}{4d}$, i.e. 
$\nu_0'd<\frac 1 4$.  Fix $T$ and define 
$$S_0=\{g\in G: \ p_g(T) \geq \nu_0\}.$$  
Since the $p_g(T)$ sum to $1$ and $\nu_0<1/d$, $S_0$ is nonempty.  
By Lemma \ref{F1}, we have that with probability at least $1- c^{-T}$, 
$$ p_h(T_1)\geq \nu_0' \quad \mbox{ for all } h \in H_0=\langle S_0 \rangle,$$ 
where $T_1=2^dT$.  
If $H_0\ne G$, then Lemma \ref{F2} applied to $H_0$ guarantees that 
$p_{G\setminus H_0}(t)\geq 1/4 > \nu_0'd$ at some time $t$; let $T_1'\geq T_1$
be the first such time.  Then there must be a $g^* \not\in H_0$ with 
$p_{g^*}(T_1')\geq \nu_0'$.
Now Lemma \ref{F3} implies that with probability at least $1-C^{-T}$, 
$p_h(T_1') \geq 2\nu_0'/3 =: \nu_1$ for all $h \in H_0$.  Let 
$$S_1=\{g\in G: \ p_g(T_1') \geq \nu_1\};$$ since $\nu_0' \geq \nu_1$, $S_1$ 
includes $g^*$ as well as all of $H_0$.  Hence $H_1=\langle S_1\rangle$
is strictly larger than $H_0$, and we may repeat the argument.  
After some number $k\leq d$ of iterations, $H_k$ must equal $G$ as desired.  

To complete the proof, note that once the $\nu_j$ have been fixed, 
the exceptional probability in this argument is 
on the order of $c^{-T}$ for some $c>1$, hence it can be made arbitrarily small 
by choosing a sufficiently large initial time $T$.
\end{proof}

\begin{proof}[Proof of Lemma \ref{F1}.]
In running the evolution (\ref{e:probrule}) until time $2T$, each 
$p_i \geq \nu /2$ for $g_i\in S$ as in 
(\ref{e:overr}).  Hence the probability at each step of adding any product 
$g_ig_j$ to the urn where $g_i,g_j\in S$ 
is at least $(\frac{\nu}{2})^2$.  By choosing $\delta=\frac{1}{2}$ in Chernoff's 
bound, the number of occurrences of $g_ig_j$ in the urn at time $2T$ exceeds 
$\frac{1}{2}T(\frac{\nu}{2})^2$ (or, equivalently, $p_{g_ig_j}(2T)\geq \frac {\nu^2} {16}$)
with probability at least $1-e^{-\frac{1}{8}T(\frac{\nu}{2})^2}$.  

Since every element of $\langle S \rangle$ can be expressed as a
product of at most $\left|\langle S \rangle \right|\leq |G|=d$ elements of
$S$, we iterate this argument $d$ times to obtain all
elements of $\langle S\rangle$.  Specifically, at time $2^dT$ each element
of $H$ has density at least
\[ 16\frac{\nu^{2^d}}{{16}^{2^d}} \]
with probability at least $1- c^{-T}$, where $c>1$ is a suitable constant depending 
on $\nu$ and $d$.
\end{proof}

\begin{proof}[Proof of Lemma \ref{F2}.]
We will couple the behavior of $n_H(t)$ and $n_{G\setminus H}(t)$ under the evolution 
(\ref{e:probrule}) with the behavior of the quantities $n_0(t)$ and $n_1(t)$ associated
to a second urn running the same evolution on the group $\Z /2\Z$, beginning with 
$n_H(T)$ balls labeled 0 and $n_{G\setminus H}(T)$ balls labeled 1.  
Picking $h\in H$ and $k\in G\setminus H$ from the first urn correspond to 
picking $0$ and $1$ from the second urn, respectively.
Since $hk$ and $kh$ are in $G\setminus H$ we 
see that at times when $n_H=n_0$ and $n_{G\setminus H}=n_1$, the probability of increasing 
$n_{G\setminus H}$ is at least as great as that of increasing $n_1$.  Coupling the urns 
at these times shows that $n_{G\setminus H}(t)\geq n_1(t)$ for all $t$.   
However, by Corollary \ref{c:1}, $n_1(t)/t$ approaches $1/2$ as $t$ increases, 
hence $n_{G\setminus H}(t)/t=p_{G\setminus H}(t)$ eventually exceeds $1/4$. 
\end{proof}

\begin{proof}[Proof of Lemma \ref{F3}.]
The reasoning is similar to that in Theorem \ref{t:1}.  We begin by showing for 
sufficiently small $\nu$ that if at time $T$, $p_h(T) \geq \nu$ for all $h\in H$ and 
$p_{G\setminus H} (t) < \nu d$ for $T\leq t < 3T/2$, then with probability
$1-C^{-T}$, $p_h(3T/2) \geq \nu$ for all $h\in H$. 
A lower bound for the probability of adding $h\in H$ to the urn at time 
$T\leq t < \frac{3T}{2}$ is
\begin{equation}\label{e:probH}
\left(p_H(t)\right)^2\left(2\frac{\nu/p_H(t)}{3/2} - d_0 \left(\frac{\nu/p_H(t)}{3/2}\right)^2\right)
\end{equation}
where the first term is the probability of 
picking both elements from $H$ and the second term is the lower bound given by 
Lemma \ref{l:1} applied with $r=3/2$ to the group $H$ 
of size $d_0$ with normalized densities at least $\nu/p_H(t)$.  
By assumption $p_H(t)$ is at least $1-d\nu$, so (\ref{e:probH}) is bounded below by
\begin{equation*}
\frac{4}{3}p_H(t) \nu - d_0 \left(\frac{2}{3}\nu\right)^2 \geq 
\frac{4\nu}{3}(1  - \frac{4}{3}\nu d)=:\mu.
\end{equation*}

Thus by Chernoff's bound, with probability at least
\[ 1- e^{-\frac{\delta ^2}{2}\mu\frac{T}{2} }, \] 
there will be at least $(1-\delta)\mu\frac{T}{2}$ balls labeled $h\in H$ added to the 
urn between times $T$ and $\frac{3T}{2}$.  For small enough $\nu$ we can choose 
$\delta$ small enough so that 
\[(1-\delta)\mu \geq \nu, \] 
and thus with probability $1- C_0^{-T}$ for some $C_0>1$, we have 
$p_h(\frac{3T}{2}) \geq \nu$ for all $h\in H$.  

By choosing $j$ so that $(\frac{3}{2})^jT \leq T_1 < (\frac{3}{2})^{j+1}T$, we can 
iterate this argument $j$ times to conclude that $p_h\left((\frac{3}{2})^jT\right)\geq \nu$ 
for all $h\in H$ with probability at least
$$ 1-\sum_{i=0}^{j-1}C_0^{-(\frac{3}{2})^iT}=1-C^{-T},$$  
where (as in the proof of Theorem \ref{t:1}) $C$ can be chosen independently
of $j$.  The argument is completed by noting that again as in (\ref{e:overr}),   
$p_h\left((\frac{3}{2})^jT\right)\geq \nu$  implies $p_h(t) \geq\frac{2}{3}\nu$ for all 
$(\frac{3}{2})^jT \leq t < (\frac{3}{2})^{j+1}T$.
\end{proof}


\vskip 0.2in
{\scriptsize
{\bf Aaron Abrams,} Department of Mathematics and Computer Science,
Emory University, Atlanta, GA  30322, and
Mathematical Sciences Research Institute, 17 Gauss Way,
Berkeley California, 94720 (abrams@msri.org)
\\
{\bf Henry Landau,}
Teachers College Columbia University (hjlandau@yahoo.com)
\\
{\bf Zeph Landau,} Department of Mathematics R8133, The City College of New York,
Convent Ave \& 138th, New York, NY  10031 (landau@sci.ccny.cuny.edu) \\ 
{\bf James Pommersheim,} Reed College, 3203 SE Woodstock Blvd.,
Portland, OR 97202 (jamie@reed.edu)\\
{\bf Eric Zaslow,} Department of Mathematics, Northwestern University,
2033 Sheridan Road, Evanston, IL  60208 (zaslow@math.northwestern.edu)
}


\end{document}